\documentclass[11pt]{article}
\usepackage{amsmath,amssymb}

\newcommand{\Tbar}{\ensuremath{\overline{\mathcal{T}}}}

\begin{document}

\title{Weil-Petersson perspectives}         
\author{Scott A. Wolpert}        
\date{March 1, 2005}          
\maketitle

\section{Introduction.}
We highlight recent progresses in the study of the Weil-Petersson (WP) geometry of finite dimensional Teichm\"{u}ller spaces.  For recent progress on and the understanding of infinite dimensional Teichm\"{u}ller spaces the reader is directed to \cite{TakTI,TakTII}.  As part of the highlight, we also present possible directions for future investigations.  Recent works on WP geometry involve new techniques which present new opportunities.

We begin with background highlights.  The reader should see \cite{Ahsome,Bersdeg,Ng, Royicm}, as well as \cite[Sec. 1 and 2]{Wlcomp} for particulars of our setup and notation for the augmented Teichm\"{u}ller space $\Tbar$, the mapping class group $MCG$, Fenchel-Nielsen (FN) coordinates and the complex of curves $C(F)$ for a 
surface $F$.  We describe the main elements.  

For a reference topological surface $F$ of genus $g$ with $n$ punctures and negative Euler characteristic, a Riemann surface $R$ (homeomorphic to $F$) with complete hyperbolic metric, a marked Riemann/hyperbolic surface is the equivalence class of a pair $\{(R,f)\}$ for $f:F\rightarrow R$ an orientation preserving homeomorphism; equivalence for post composition with a conformal homeomorphism.  Teichm\"{u}ller space $\mathcal T$ is the space of equivalence classes, \cite{Ahsome, ImTan, Ng}.  The mapping class group $MCG=Homeo^+(F)\slash Homeo_0(F)$ acts properly discontinuously on $\mathcal T$ by taking $\{(R,f)\}$ to $\{(R,f\circ h^{-1})\}$ for a homeomorphism $h$ of $F$.  

$\mathcal T$ is a complex manifold with the $MCG$ acting by biholomorphisms.  The cotangent space at $\{(R,f)\}$ is $Q(R)$, the space of integrable holomorphic quadratic differentials, \cite{Ahsome, Ng}.  For $\phi,\psi\in Q(R)$ the Teichm\"{u}ller (Finsler metric) conorm of $\phi$ is $\int_R|\phi|$, while the WP dual Hermitian pairing is $\int_R\phi\overline{\psi}(ds^2)^{-1}$, where $ds^2$ is the $R$ complete hyperbolic metric.  The WP metric is K\"{a}hler, not complete, with negative sectional curvature and the $MCG$ acts by isometries (see Section 3 below), \cite{Ahsome, Msext,  Royicm, Trcurv, Wlchern}. The {\em complex of curves} $C(F)$ is defined as follows.  The vertices of $C(F)$ are free homotopy classes of homotopically nontrivial, nonperipheral, simple closed curves on $F$.  A $k$-simplex of $C(F)$ consists of $k+1$ distinct homotopy classes of mutually disjoint simple closed curves.  A maximal set of mutually disjoint simple closed curves, a {\em pants decomposition}, has $3g-3+n$ elements for $F$ a genus $g$, $n$ punctured surface.  

A free homotopy class $\gamma$ on $F$ determines the geodesic-length function $\ell_{\gamma}$ on $\mathcal T$; $\ell_{\gamma}(\{R,f,)\})$ is defined to be the hyperbolic length of the unique geodesic freely homotopic to $f(\gamma)$.  The geodesic-length functions are convex along WP geodesics, \cite{Wlnielsen, Wlreprise}.   On a hyperbolic surface the geodesics for a pants decomposition decompose the surface into geometric pairs of pants: subsurfaces homeomorphic to spheres with a combination of three discs or points removed, \cite{Abbook, ImTan}.  The hyperbolic geometric structure on a pair of pants is uniquely determined by its boundary geodesic-lengths.  Two pants boundaries with a common length can be abutted to form a new hyperbolic surface (a complete hyperbolic structure with possible further geodesic boundaries.)  The common length $\ell$ for the joined boundary and the offset, or {\em twist} $\tau$, for adjoining the boundaries combine to provide parameters $(\ell,\tau)$ for the construction.  The twist $\tau$ is measured as displacement along boundaries.  The Fenchel-Nielsen (FN) parameters $(\ell_{\gamma_1},\tau_{\gamma_1},\dots,\ell_{\gamma_{3g-3+n}},\tau_{\gamma_{3g-3+n}})$ valued in $(\mathbb R_+\times\mathbb R)^{3g-3+n}$ provide global real-analytic coordinates for $\mathcal T$, \cite{Abbook, MasFN, WlFN}.  Each pants decomposition determines a global coordinate.  

A bordification of $\mathcal T$, the {\em augmented Teichm\"{u}ller space}, is introduced by extending the range of the parameters. For an $\ell_{\gamma}$ equal to zero, the twist is not defined and in place of the geodesic for $\gamma$ there appears a pair of cusps.  Following Bers \cite{Bersdeg} the extended FN parameters describe marked (possibly) noded Riemann surfaces (marked stable curves.)  An equivalence relation is defined for marked noded Riemann surfaces and a construction is provided for adjoining to $\mathcal T$ frontier spaces (where subsets of geodesic-lengths vanish) to obtain the augmented Teichm\"{u}ller space $\Tbar$, \cite{Abdegn, Abbook}. $\Tbar$ is not locally compact since in a neighborhood of $\ell_{\gamma}$ vanishing the FN angle $\theta_{\gamma}=2\pi\tau_{\gamma}/\ell_{\gamma}$ has values filling $\mathbb R$.  The leading-term expansion for the WP metric is provided in \cite{DW2, Yam, Wlcomp}. Following Masur the WP metric extends to a complete metric on $\Tbar$, \cite{Msext}.  

The group $MCG$ acts (not properly discontinuously) as a group of homeomorphisms of $\Tbar$ and $\Tbar\slash MCG$ is topologically the Deligne-Mumford compactified moduli space of stable curves $\overline{\mathcal M}$.  We include the empty set as a $-1$-simplex of $C(F)$.  The complex is partially ordered by inclusion of simplices.  For a marked noded Riemann surface $\Lambda(\{(R,f)\})\in C(F)$ is the simplex of free homotopy classes on $F$ mapped to nodes of $R$.  The level sets of $\Lambda$ are the strata of $\Tbar$.  

$(\Tbar,d_{WP})$ is a $CAT(0)$ metric space; see \cite{DW2, Yam, Wlcomp} and the attribution to Benson Farb in \cite{MW}.  The general structure of $CAT(0)$ spaces is described in detail in \cite{BH}.  $CAT(0)$ spaces are complete metric spaces of {\em nonpositive curvature.}   In particular for $p,q\in\Tbar$ there is a unique length-minimizing path $\widehat{pq}$ connecting $p$ and $q$.  An additional property is that WP geodesics do not 
{\em refract}; an open WP geodesic segment $\widehat{pq}-\{p,q\}$ is contained in the (open) stratum $\Lambda(p)\cap\Lambda(q)$ \cite{DW2,Yam, Wlcomp}; an open segment is a solution of the WP geodesic differential equation on a product of Teichm\"{u}ller spaces.  The stratum of $\Tbar$ (the level sets of $\Lambda$) are totally geodesic complete subspaces.  Strata also have a metric-intrinsic description: a stratum is the union of all open length-minimizing paths containing a given point, \cite[Thrm. 13]{Wlcomp}. 

A classification of {\em flats}, the locally Euclidean isometrically embedded subspaces, of $\Tbar$ is given in \cite[Prop. 16]{Wlcomp}.  Each flat is contained in a proper substratum and is the Cartesian product of geodesics from component Teichm\"{u}ller spaces.  A structure theorem characterizes limits of WP geodesics on $\Tbar$ modulo the action of the mapping class group, \cite[Sec. 7]{Wlcomp}.  Modulo the action the general limit is the unique length-minimizing piecewise geodesic path connecting the initial point, a sequence of strata and the final point.  An application is the construction of  axes in $\Tbar$ for elements of the $MCG$, \cite{DW2,Wlcomp}.  The geodesic limit behavior is suggested by studying the sequence $\widehat{pT_{\gamma}^np}$ for a Dehn twist $T_{\gamma}\in MCG$: modulo $MCG$ the limit is two copies of the length-minimizing path from $p$ to $\{\ell_{\gamma}=0\}$. 

Jeffrey Brock established a collection of important results on WP synthetic geometry, \cite{Brkwp, Brkwpvs}. Brock introduced an approach for approximating WP geodesics.  A result is that the geodesic rays from a point of $\mathcal T$ to the maximally noded Riemann surfaces have initial tangents dense in the initial tangent space.  We used the approach in \cite[Coro. 19]{Wlcomp} to find that $\Tbar$ is the closed convex hull of the discrete subset of marked maximally noded Riemann surfaces.  We also showed that the geodesics connecting the maximally noded Riemann surfaces have tangents dense in the tangent bundle of $\mathcal T$. 

Howard Masur and Michael Wolf theorem established the WP counterpart to H. Royden's celebrated theorem: each WP isometry of $\mathcal T$ is induced by an element of the extended $MCG$, \cite{MW}.  
A simplified proof of the Masur-Wolf result is given as follows: \cite[Thrm. 20]{Wlcomp}, the elements of $Isom_{WP}(\mathcal T)$ extend to isometries of $\Tbar$; the extensions preserve the metric-intrinsic stratum of $\Tbar$ and so correspond to simplicial automorphisms of $C(F)$; the simplicial automorphisms of $C(F)$ are in general induced by the elements of the $MCG$ from the work of N. Ivanov, M. Korkmaz and F. Luo, \cite{Ivaut,Krk, Luaut}; and finally $\Tbar-\mathcal T$ is a uniqueness set for WP isometries.  Brock and Dan Margalit have recently extended available techniques to include the special $(g,n)$ types of $(1,1), (1,2)$ and $(2,0)$, \cite{BrMr}.

I would like to take the opportunity to thank Benson Farb, Richard Wentworth and Howard Weiss for their assistance.
 
\section{Classical metrics for $\mathcal T$.}
A puzzle of Teichm\"{u}ller space geometry is that the classical Teichm\"{u}ller and Weil-Petersson metrics in combination have desired properties: K\"{a}hler, complete, finite volume, negative curvature, and a suitable sphere at infinity, yet individually the metrics already lack the basic properties of completeness and non-positive curvature.  In fact N. Ivanov showed that $\mathcal M_{g,n}$, $3g-3+n\ge 2$, admits no complete Riemannian metric of pinched negative sectional curvature, \cite{Ivcurv}.  Recently authors have considered asymptotics of the WP metric and curvature with application to introducing modifications to obtain {\em designer metrics}.  

Curt McMullen considered adding a small multiple of the Hessian of the {\em short geodesic log length sum} to the WP K\"{a}hler form
$$
\omega_{1/\ell}=\omega_{WP}+i\,c\sum_{\ell_{\gamma}<\epsilon}\partial\overline{\partial}Log\,\ell_{\gamma}
$$
(for $Log \approx \min{\{\log, 0}\}$) to obtain the K\"{a}hler form for a modified metric , \cite{McM}. The constructed metric is {\em K\"{a}hler hyperbolic}: on the universal cover the K\"{a}hler form $\omega_{1/\ell}$ is the exterior derivative of a bounded $1$-form and the injectivity radius is positive; on the moduli space $\mathcal M_{g,n}$ the metric $g_{1/\ell}$ is: complete, finite volume, of bounded sectional curvature.  Of primary interest is that the constructed metric $g_{1/\ell}$ is comparable to the Teichm\"{u}ller metric.  The result provides that the Teichm\"{u}ller metric qualitatively has the properties of a K\"{a}hler hyperbolic metric, \cite{McM}. As applications McMullen establishes a positive lower bound for the Teichm\"{u}ller Rayleigh-Ritz quotient and  a 
Teichm\"{u}ller metric isoperimetric inequality for complex manifolds.  McMullen also uses the metric to give a simple derivation that the sign of the orbifold Euler characteristic alternates with the parity of the complex dimension. The techniques of \cite{Wlreprise} can be applied to show that in a neighborhood of the compactification divisor 
$\mathcal D=\overline{\mathcal M_g}-\mathcal M_g\subset\overline{\mathcal M_g}$  the McMullen modification to WP is principally in the directions transverse to the divisor.  

Lipman Bers promoted the question of understanding the equivalence of the classical metrics on Teichm\"{u}ller spaces in 1972, \cite{Berssurv}.  A collection of authors have now achieved major breakthroughs.    Kefeng Liu, Xiaofeng Sun and Shing-Tung Yau have analyzed WP curvature, applied Yau's Schwarz lemma, and considered the Bers embedding to compare and study the K\"{a}hler-Einstein metric for Teichm\"{u}ller space and also the classical metrics, \cite{LSY}. Sai-Kee Yeung considered the Bers embedding and applied Schwarz type lemmas to compare classical metrics, \cite{Yun}.  Yeung's approach is based on his analysis of convexity properties of fractional powers of geodesic-length functions, \cite{Yun1}.
Bo-Yong Chen considered McMullen's metric, the WP K\"{a}hler potential of Takhtajan-Teo,
and $L^2$-estimates to establish equivalence of the Bergman and Teichm\"{u}ller metrics, \cite{Chn}.

In earlier work, Cheng-Yau \cite{ChYa} and Mok-Yau \cite{MoYa} had established existence and
completeness for a K\"{a}hler-Einstein metric for Teichm\"{u}ller space.  Liu-Sun-Yau study the negative WP Ricci form $-Ricci_{WP}$; the form defines a complete K\"{a}hler metric on the Teichm\"{u}ller space.  The authors find that the metric is equivalent to the {\em asymptotic Poincar\'{e} metric} and that asymptotically its holomorphic sectional 
curvature has negative upper and lower bounds.  The authors develop the asymptotics for the curvature of the metric.  

To overcome difficulty with controlling interior curvatures Liu-Sun-Yau further introduce the modification 
$$
g_{LSY}=-Ricci_{WP}+c\,g_{WP},\  c>0,\  \cite{LSY}.
$$
The authors find that the LSY metric is complete with bounded negative holomorphic sectional curvature and bounded negative Ricci curvature.  The authors apply their techniques and results of others to establish comparability of the seven classical metrics
\begin{quote}
LSY $\sim$ asymptotic Poincar\'{e}   $\sim$ K\"{a}hler-Einstein $\sim$ Teichm\"{u}ller-Kobayashi $\sim$ McMullen 
 $\sim$ Bergman $\sim$ Carath\'{e}odory.
\end{quote}
Yeung established the comparability of the last five listed classical metrics.  The classical metrics are now expected to have the same qualitative behavior.  In combination the metrics have a substantial list of properties.  Liu-Sun-Yau use their understanding of the metrics in a neighborhood of the compactification divisor $\mathcal D\subset\overline{\mathcal M_g}$ to show that the logarithmic cotangent bundle (Bers' bundle of {\em holomorphic $2$-differentials}) of the compactified moduli space is stable in the sense of Mumford.  The authors also show that the bundle is stable with respect to its first Chern class.  The authors further find that the K\"{a}hler-Einstein metric has bounded geometry in a strong sense.  The comparison of metrics now provides a major opportunity to combine different approaches for studying Teichm\"{u}ller geometry. 

\section{WP synthetic geometry.}
A collection of authors including Jeffrey Brock, Sumio Yamada and Richard Wentworth have recently studied WP synthetic geometry.  A range of new techniques have been developed.
Brock established a collection of very interesting results on the
large-scale behavior of WP distance.  Brock considered the {\em pants graph}
$C_{\mathbf P}(F)\subset C(F)$ having vertices the distinct pants decompositions of $F$
and joining edges of unit-length for pants decompositions differing by an
elementary move,\cite{Brkwp}.  He showed that the $0$-skeleton of $C_{\mathbf
P}(F)$ with the edge-metric is quasi-isometric to  $\mathcal T$ with the WP metric.  He further
showed for $p,q\in \mathcal T$ that the corresponding quasifuchsian hyperbolic
three-manifold has convex core volume comparable to $d_{WP}(p,q)$.  At
large-scale, WP distance and convex core volume are approximately
combinatorially determined.  

Brock's approach begins with Bers' observation that there is a constant $L$ depending only on the 
topological type of $F$, such that a surface $R$ has a pants decomposition with lengths $\ell_{\gamma}(R)<L$.  For $\mathcal P$ a pants decomposition Brock associates the sublevel set 
$$
V(\mathcal P)=\{R\,|\,\gamma\in \mathcal P,\ell_{\gamma}<L\}.
$$  
By Bers' theorem the sublevel sets $V(\mathcal P)$ cover $\mathcal T$, \cite{Bersdeg}.  Brock shows that the configuration of the sublevel sets in terms of WP distance on $\mathcal T$ is coarsely approximated by the metric space $C_{\mathbf P}(F)$.  He shows for $R,S,\in\mathcal T$ with $R\in V(\mathcal P_R),\,S\in V(\mathcal P_S)$ that 
$$
d_{\mathbf P}(\mathcal P_R,\mathcal P_S)\asymp d_{WP}(R,S)
$$ 
for $d_{\mathbf P}$ the $C_{\mathbf P}(F)$ edge-metric, \cite{Brkwp}. Brock's underlying idea is that $C_{\mathbf P}(F)$ is the $1$-skeleton for the nerve of the covering $\{V(\mathcal P)\}$ and that the minimal count of {\em elementary moves} in $C_{\mathbf P}(F)$ approximates WP distance.   The paths in $C_{\mathbf P}(F)$ can be encoded as sequences of elementary moves.  Can sequences of elementary moves be used to study the bi-infinite WP geodesics?  What is the large $g,\,n$ behavior of $diam_{WP}(\mathcal M_{g,n})$?

Questions regarding WP geodesics include behavior in-the-small, as well as in-the-large.  Understanding the WP {\em metric tangent cones} of $\Tbar$, as well as the {\em metric tangent cone bundle} over $\Tbar$ is a basic matter.  Recall the definition of the Alexandrov angle for unit-speed geodesics $\alpha(t)$ and $\beta(t)$ emanating from a common point 
$\alpha(0)=\beta(0)=o$ of a metric space $(X,d)$, \cite[pg. 184]{BH}.  The Alexandrov angle between 
$\alpha$ and $\beta$ is defined by 
$$
\angle(\alpha,\beta)=\lim_{t\rightarrow 0}\frac{1}{2t}d(\alpha(t),\beta(t)).
$$  
Basic properties of the angle are provided in the opening and closing sections of \cite[Chap. II.3]{BH}.  
Angle zero provides an equivalence relation on the space of germs of unit-speed geodesics beginning at $o$.  For the model space 
$(\mathbb H;\,ds^2=dr^2+r^6d\theta^2, r>0)$ for the WP metric transverse to a stratum: the geodesics $\{\theta=\mbox{constant}\}$ are at angle zero at the {\em origin point} $\{r=0\}$. The WP metric tangent cone at $p\in\Tbar$ is defined as the space of germs of constant-speed geodesics beginning at $p$ modulo the relation of 
Alexandrov angle zero.  Part of the investigation would be to define the metric tangent cone bundle.  
Question: can WP geodesics spiral to a stratum?  In particular for $\mathcal S$ a stratum defined by the vanishing of a geodesic-length function $\ell_{\gamma}$ and $\theta_{\gamma}$ the corresponding Fenchel-Nielsen angle, is $\theta_{\gamma}$ bounded along each geodesic ending on $\mathcal S$ (see the next section for results of Wentworth precluding spiraling for harmonic maps from $2$-dimensional domains.)

There are questions regarding in-the-large behavior of WP geodesics.  A basic invariant of hyperbolic metrics, the relative systole $sys_{rel}(R)$ of a surface $R$ is the length of the shortest closed geodesic.  In \cite{Wlcomp} we note that the WP injectivity radius $inj_{WP}$ of the moduli space $\mathcal M$ is comparable to $(sys_{rel})^{1/2}$.  A basic question is to understand the relative systole $sys_{rel}$ along infinite WP geodesics in $\mathcal T$.  It would be interesting to understand the combinatorics of the sequence of short geodesics, and any limits of the (weighted) collections of simple closed geodesics for which the relative systole is realized.  Brock, Masur and Minsky have preliminary results on such limits in the space of measured geodesic laminations $\mathcal{MGL}$.  The authors are studying the family of hyperbolic metrics $\{R_t\}$ over a WP geodesic.  A geometric model $M$ is constructed from the family; the authors seek to show that the model is biLipschitz to the {\em thick part} of a hyperbolic $3$-dimensional structure.  To construct $M$, the authors remove the {\em thin parts} of the fibers and then use the hyperbolic metric along fibers and (subsurface-component rescalings) of the WP metric in the parameter to define a $3$-dimensional metric.   

A second matter is to analyze the rate of convergence of pairs consisting of an infinite WP geodesic and either a second infinite WP geodesic or a stratum.  Although the WP metric on $\mathcal T$ has negative sectional curvature there are directions of almost zero curvature at the bordification $\Tbar-\mathcal T$, as well as flats in the bordification.  To clarify, Zheng Huang has shown that the sectional curvature on $\mathcal T$ is bounded as $-c'\,(sys_{rel})^{-1}<sec_{WP}<-c''\, sys_{rel}$ for positive constants, \cite{Zh2}.  The divergence of a geodesic and a geodesic/stratum is a measure of the curvature of the join.  

A third matter proposed by Jeffrey Brock is to investigate the dynamics of the WP geodesic flow on the tangent bundle $T\mathcal T$.  The $CAT(0)$ geometry can be used to show that the {\em finite and semi-finite} WP geodesics form a subset of smaller Hausdorff dimension, \cite[Sec. 6]{Wlcomp}.   Are there WP geodesics with lifts dense in $T\mathcal T$?  Are the lifts of  axes of (pseudo Anosov) elements of the $MCG$ dense in $T\mathcal T$?  What is the growth rate of the translation lengths for the conjugacy classes of the pseudo Anosov elements?    

A more general matter is to understand the behavior of quasi-geodesics and especially quasi-flats, quasi-isometric embeddings of Euclidean space into $\mathcal T$.  Brock and Farb have conjectured that the maximal dimension of a quasi-flat is $\lfloor\frac{g+n}{2}\rfloor-1$ (the quantity does provide a lower bound \cite{BF, Wlcomp}.)  The maximal dimension of a quasi-flat is the rank in the sense of Gromov.  The rank is important for understanding the global WP geometry, as well as for understanding mapping class groups.  Brock-Farb, Jason Behrstock and Javier Aramayona have found that certain low dimensional Teichm\"{u}ller spaces are Gromov-hyperbolic and necessarily rank one, \cite{ Arm,Bhr,BF}.

\section{Harmonic maps to $\Tbar$.}
A collection of authors including Georgios Daskalopoulos, Richard Wentworth and Sumio Yamada have considered harmonic maps from (finite volume) Riemannian domains $\Omega$ into the $CAT(0)$ space $\Tbar$, \cite{DW1, DW2, Wreg, Yam}.  The authors apply the Sobolev theory of Korevaar-Schoen, \cite{KShar}, to study maps energy-minimizing for prescribed boundary values. Possible applications are rigidity results for homomorphisms of lattices in Lie groups to mapping class groups and existence results for harmonic representatives of the classifying maps associated to symplectic Lefschetz pencils, \cite{DW2}.  

Understanding the behavior of a harmonic map to a neighborhood of a lower stratum of $\Tbar$ is a basic question.  Again this is especially important since $\Tbar$ is not locally compact in a neighborhood of a lower stratum; in each neighborhood the corresponding Fenchel-Nielson angles surject to $\mathbb R$.  A harmonic map or even a WP geodesic may {\em spiral} with unbounded Fenchel-Nielsen angles.  
With the non-refractions of WP geodesics at lower stratum \cite{DW2, Yam, Wlcomp}, a question is to consider analogs of {\em non-refraction} for harmonic mappings.  In particular for the complex of curves $C(F)$ and $\Lambda:\Tbar\rightarrow C(F)$ the labeling function and $u:\Omega\rightarrow\Tbar$ a harmonic map, what is the behavior of the composition $\Lambda\circ u$? For $2$-dimensional domains Wentworth has already provided important results on {\em non-spiraling} and on $\Lambda \circ u$ being constant on the interior of $\Omega$, \cite{Wreg}.   A second general matter is to develop an encompassing approach (including treating regularity and singularity behavior) for approximation of maps harmonic to a neighborhood of a stratum by harmonic maps to model spaces, in particular for: $(\mathbb H;\,ds^2=dr^2+r^6d\theta^2, r>0$).

\section{Characteristic classes and the WP K\"{a}hler form.}
Andrew McIntyre in joint work with Leon Takhtajan \cite{McTk}, Maryam Mirzakhani \cite{Mirgrow, Mirvol, Mirwitt}, and Lin Weng in part in joint work with Wing-Keung To \cite{Weng} have been studying the questions in algebraic geometry involving the WP K\"{a}hler form.  The authors' considerations are guided by the study of the Quillen and Arakelov metrics, and in particular calculations from {\em string theory.}  We sketch aspects of their work and use the opportunity to describe research with Kunio Obitsu.

Andrew McIntyre and Leon Takhtajan have extended the work of Peter Zograf and provided a new {\em holomorphic factorization formula} for the regularized determinant $\det'\Delta_k$ of the hyperbolic Laplace operator acting on (smooth) symmetric-tensor $k$-differentials for compact Riemann surfaces, \cite{McTk}.  Alexey Kokotov and Dmitry Korotkin have also provided a new {\em holomorphic factorization formula} for the regularized determinant $\det'\Delta_0$ of the hyperbolic Laplace operator acting on functions, \cite{KKo}.  Each factorization involves the exponential of an {\em action integral}: in the first formula for a Schottky uniformization, and in the second formula for a branched covering of $\mathbb CP^1$. 
The McIntryre-Takhtajan formula provides a {\em $\partial\overline{\partial}$ antiderivative} for the celebrated families local index theorem
$$
\overline{\partial}\partial\log\frac{\det N_k}{\det'\Delta_k}=\frac{6k^2-6k+1}{6\pi i}\ \omega_{WP}
$$
for $N_k$ the Gram matrix of the natural basis for holomorphic $k$-differentials relative to the Petersson product, \cite{McTk}.  The McIntryre-Takhtajan formula for $\frac{\det N_k}{\det'\Delta_k}$ gives rise to an isometry between the determinant bundle for holomorphic $k$-differentials with Quillen metric and with a metric defined from the Liouville action.  Recall that the quotient $\frac{\det N_k}{\det'\Delta_k}$ is the Quillen norm of the natural frame for $\det N_k$.  The formulas represent progress in the ongoing study of the behavior of the Quillen metric, the hyperbolic regularized determinant $\det'\Delta_k$ and positive integral values of the Selberg zeta function, since $\det'\Delta_k=c_{g,k}Z(k), k>2$, \cite{DPh, Sardet}.  The above formulas provide another approach for studying the degeneration of $\det'\Delta_k$; see \cite{Hejreg,Wlselb}.

Mumford's tautological class $\kappa_1$, the pushdown of the square of the relative dualizing sheaf from the universal curve, is represented by the WP class $\frac{1}{2\pi^2}\,\omega_{WP}$, \cite{Wlhyp}.  The top  self-intersection number of $\kappa_1$ on $\overline{\mathcal M_{g,n}}$ is the WP volume.  From effective estimates for intersections of divisors on the moduli space, Georg Schumacher and Stefano Trapani \cite{SchTr} have given lower bounds for $vol_{WP}(\overline{\mathcal M_{g,n}})$.  From Robert Penner's \cite{Pen} decorated Teichm\"{u}ller theory and a combinatorial description of the moduli space, Samuel Grushevsky \cite{Gru} has given upper bounds with the same leading growth order.

In a series of innovative papers Maryam Mirzakhani has presented a collection of new results on hyperbolic geometry and calculations of WP integrals.  The asymptotics for the count of the number of simple closed geodesics on a hyperbolic surface of at most a given length is presented in \cite{Mirgrow}.  She establishes the asymptotic
$$
\#\{\gamma\,|\,\ell_{\gamma}(R)\le L\}\sim c_R\, L^{6g-6+2n}.
$$    
A recursive method for calculating the WP volumes of moduli spaces of bordered hyperbolic surfaces with prescribed boundary lengths is presented in \cite{Mirvol}. And a proof of the Witten-Kontsevich formula for the tautological classes on $\overline{\mathcal M_{g,n}}$ is presented in \cite{Mirwitt}.  

Central to Mirzakhani's considerations is a recursive scheme for evaluating WP integrals over the moduli space of bordered hyperbolic surfaces with prescribed boundary lengths.  Her approach is based on recognizing an integration role for McShane's length sum identity.  
Gregg McShane discovered a remarkable identity for geodesic-lengths of simple closed curves on punctured hyperbolic surfaces, \cite{McSh}.  To illustrate the recursive scheme for evaluation of integrals we sketch the consideration 
for $(g,n)=(1,1)$.  

For the $(1,1)$ case the identity provides that 
$$
\sum_{\gamma\ scg}\frac{1}{1+e^{\ell_{\gamma}}}=\frac12
$$ 
for the sum over simple closed geodesics (scg's). The identity corresponds to a decomposition of a horocycle about the puncture; the identity arises from classifying the behavior of simple complete geodesics emanating from the puncture.   Mirzakhani's insight is that the identity can be combined with the $d\tau\wedge d\ell$ formula for $\omega_{WP}$, \cite{Wldtau} ,to give an {\em unfolding} of the $\mathcal M_{1,1}$ volume integral.  For a $MCG$-fundamental domain $\mathcal F_{1,1}\subset\mathcal T_{1,1}$, she observes that
\begin{multline*}
\int_{\mathcal M_{1,1}}\frac12\,\omega_{WP}=\int_{\mathcal F_{1,1}}\sum_{h\in MCG\slash Stab_{\gamma}}\frac{1}{1+e^{\ell_{h(\gamma)}}}\,d\tau\wedge d\ell =\\
\sum_{h\in MCG\slash Stab_{\gamma}}\int_{h^{-1}(\mathcal F_{1,1})}\frac{1}{1+e^{\ell}}\,d\tau\wedge d\ell = \int_{\mathcal T_{1,1}\slash Stab_{\gamma}}\frac{1}{1+e^{\ell}}\, d\tau\wedge d\ell
\end{multline*}
(using that $\ell_{\gamma}\circ h^{-1}=\ell_{h(\gamma)}$) with the last integral elementary since 
$$
\mathcal T_{1,1}\slash Stab_{\gamma}=\{\ell>0,0<\tau<\ell\}.
$$

Mirzakhani established a general identity for bordered hyperbolic surfaces that generalizes McShane's identity, \cite[Sec. 4]{Mirvol}.   The general identity is based on a sum over configurations of simple closed curves ({\em sub-partitions}) which with a fixed boundary bound a pair of pants.  Mirzakhani uses the identity to {\em unfold} WP integrals to sums of product integrals over lower dimensional moduli spaces for surfaces with boundaries.  Her overall approach provides a recursive scheme for determining WP volumes $vol_{WP}(\mathcal M_g(b_1,\dots,b_n))$ for the (real analytic) moduli spaces of hyperbolic surfaces with prescribed boundary lengths $(b_1,\dots,b_n)$.  As an instance it is shown that 
$$
vol_{WP}(\mathcal M_1(b))=\frac{\pi^2}{6}+\frac{b^2}{24}
$$ 
and 
$$
vol_{WP}(\mathcal M_1(b_1,b_2))=\frac{1}{384}(4\pi^2+b_1^2+b_2^2)(12\pi^2+b_1^2+b_2^2).
$$  
In complete generality Mirzakhani found that the WP volume is a polynomial 
$$
vol_g(b)=\sum_{|\alpha|\le 3g-3+n} c_{\alpha}\, b^{2\alpha},\ c_{\alpha}>0, \ c_{\alpha}\in\pi^{6g-6+2n-2}\mathbb Q,
$$ 
for $b$ the vector of boundary lengths and $\alpha$ a multi index, \cite{Mirvol}.  An easy application is an expansion for the volume of the tube $\mathcal N_{\epsilon}(\mathcal D)\subset \overline{\mathcal M_g}$ about the compactification divisor.  Recently Ser Peow Tan, Yan Loi Wong and Ying Zhang have also generalized McShane's identity for conical hyperbolic surfaces, \cite{TanWZ}.  A general identity is obtained by studying gaps formed by simple normal
geodesics emanating from a distinguished cone point, cusp or boundary geodesic. 

Mirzakhani also established that the volumes $vol_{WP}(\mathcal M_g(b_1,\dots,b_n))$ are determined from the intersection numbers of tautological characteristic classes on $\overline{\mathcal M_{g,n}}$.  A point of $\overline{\mathcal M_{g,n}}$ describes a Riemann surface $R$ possibly with nodes with distinct points $x_1,\dots,x_n$.  The line bundle $\mathcal L_i$ on $\overline{\mathcal M_{g,n}}$ is the unique line bundle whose fiber over $(R;x_1,\dots,x_n)$ is the cotangent space of $R$ at $x_i$; write $\psi_i=c_1(\mathcal L_i)$ for the Chern class.  Mirzakhani showed using symplectic reduction (for an $S^1$ quasi-free action following Guillemin-Sternberg) that
$$
vol_g(b)=\sum_{|\alpha|\le N}\frac{b_1^{2\alpha_1}\cdots b_n^{2\alpha_n}}{2^{|\alpha|}\alpha!(N-|\alpha|)!}\int_{\overline{\mathcal M_{g,n}}}\psi^{\alpha_1}\cdots\psi_n^{\alpha_n}\,\omega_{WP}^{N-|\alpha|}
$$
for $N=\dim_{\mathbb C}\overline{\mathcal M_{g,n}}$, \cite{Mirvol}.  She then combined the above formula and her recursive integration scheme to find that the collection of integrals
$$
\bigl<\tau_{k_1}\cdots\tau_{k_n}\bigr> = \int_{\overline{\mathcal M_{g,n}}}\psi_1^{k_1}\cdots\psi_n^{k_n}, \ \sum k_i=\dim_{\mathbb C}\overline{\mathcal M_{g,n}},
$$
satisfy the recursion for the {\em string equation} and the {\em dilaton equation}, \cite[Sec. 6]{Mirwitt}.  The recursion is the Witten-Kontsevich conjecture.  Mirzakhani's results represent major progress for the study of volumes and intersection numbers; her results clearly raise the prospect of further insights.  Can Mirzakhani's approach be applied for other integrals?  Can the considerations of Grushevsky, Schumacher-Trapani be extended to give further effective intersection estimates?

More than a decade ago Leon Takhtajan and Peter Zograf studied the local index theorem for a family of $\overline{\partial}$-operators on Riemann surfaces of type $(g,n)$, \cite{TkZg}. The authors calculated the first Chern form of the determinant line bundle provided with Quillen's metric.  For a Riemann surface with punctures there are several candidates for $\det'\Delta$ (a renormalization is necessary since punctures give rise to {\em continuous spectrum} for $\Delta$).  The authors considered the Selberg zeta function $Z(s)$ and set $\det'\Delta_k=c_{g,n}\,Z(k), k\ge2,$ \cite{TkZg}.  For $\lambda_k=\det N_k$, the determinant line bundle of the bundle of holomorphic $k$-differentials (with poles allowed at punctures), the authors found for the Chern form
$$
c_1(\lambda_k)=\frac{6k^2-6k+1}{12\pi^2}\,\omega_{WP}-\frac19\,\omega_*
$$
with the $2$-form $\omega_*$ on holomorphic quadratic differentials $\phi,\psi\in Q(R)$
$$
\omega_*(\phi,\psi)=\sum_{i=1}^n\Im\int_R\phi\overline{\psi}E_i(z;2)\,(ds^2_{hyp})^{-1}
$$
with $E_i(z;2)$, the Eisenstein-Maass series at $s=2$ for the cusp $x_i$.  By construction $\omega_*$ is a closed $(1,1)$ form; the associated Hermitian pairing (absent the imaginary part) is a $MCG$-invariant K\"{a}hler metric 
for $\mathcal T$.   The Takhtajan-Zograf (TZ) metric $g_{TZ}$ has been studied in several works,  
\cite{Obit1, Obit2, Weng}.  Kunio Obitsu showed that the metric is not complete \cite{Obit2} and is now studying the degeneration of the metric.  His estimates provide a local comparison for the TZ and WP metrics. Although beyond our exposition, we cite the important work of Lin Weng; he has been pursuing an Arakelov theory for punctured Riemann surfaces and has obtained a global comparison for the WP and TZ metrized determinant line bundles $\Delta_{WP}^{\otimes\, n^2}\le \Delta_{TZ}^{\otimes((2g-2+n)^2)}$, \cite{Weng}.  

In joint work with Obitsu we are considering the expansion for the {\em tangential to the compactification divisor} $\mathcal D\subset\overline{\mathcal M_g}$ component $g_{WP}^{tgt}$ of the WP metric, \cite{WlO}.  The tangential component is the (orthogonal) complement to Yamada's normal form $dr^2+r^6d\theta^2$ for the transversal component.  The tangential component expansion is given for a neighborhood of $\mathcal D\subset\overline{\mathcal M_g}$ for $g_{WP}$ restricted to subspaces parallel to $\mathcal D$.  For a family $\{R_{\ell}\}$ of hyperbolic surfaces given by {\em pinching} short geodesics all with common length $\ell$, we find
$$
g_{WP}^{tgt}(\ell)=g_{WP}^{tgt}(0)\,+\,\frac{\ell^2}{3}g_{TZ}(0)\,+\,O(\ell^3).
$$
The formula establishes a direct relationship between the WP and TZ metrics.  There is also a relationship with the work of Mirzakhani.  The formula is based on an explicit form of the earlier 2-term expansion for the degeneration of hyperbolic metrics, \cite[Exp. 4.2]{Wlhyp}.  What is the consequence of the expansion for the determinant line bundle?  The general question is to explore the properties of the TZ metric and its relationship to WP geometry.

\bibliographystyle{alpha}


\begin{thebibliography}{DKW00}

\bibitem[Abi77]{Abdegn}
William Abikoff.
\newblock Degenerating families of {R}iemann surfaces.
\newblock {\em Ann. of Math. (2)}, 105(1):29--44, 1977.

\bibitem[Abi80]{Abbook}
William Abikoff.
\newblock {\em The real analytic theory of {T}eichm\"uller space}.
\newblock Springer, Berlin, 1980.

\bibitem[Ahl61]{Ahsome}
Lars~V. Ahlfors.
\newblock Some remarks on {T}eichm\"uller's space of {R}iemann surfaces.
\newblock {\em Ann. of Math. (2)}, 74:171--191, 1961.

\bibitem[Ara04]{Arm}
Javier Aramayona.
\newblock The {W}eil-{P}etersson geometry of the five-times punctured sphere.
\newblock preprint, 2004.

\bibitem[Beh]{Bhr}
Jason~A. Behrstock.
\newblock Asymptotic geometry of of the mapping class group and {T}eichm\"uller
  space.
\newblock in preparation.

\bibitem[Ber72]{Berssurv}
Lipman Bers.
\newblock Uniformization, moduli, and {K}leinian groups.
\newblock {\em Bull. London Math. Soc.}, 4:257--300, 1972.

\bibitem[Ber74]{Bersdeg}
Lipman Bers.
\newblock Spaces of degenerating {R}iemann surfaces.
\newblock In {\em Discontinuous groups and Riemann surfaces (Proc. Conf., Univ.
  Maryland, College Park, Md., 1973)}, pages 43--55. Ann. of Math. Studies, No.
  79. Princeton Univ. Press, Princeton, N.J., 1974.

\bibitem[BF01]{BF}
Jeffrey~F. Brock and Benson Farb.
\newblock Curvature and rank of {T}eichm\"{u}ller space.
\newblock preprint, 2001.

\bibitem[BH99]{BH}
Martin~R. Bridson and Andr{\'e} Haefliger.
\newblock {\em Metric spaces of non-positive curvature}.
\newblock Springer-Verlag, Berlin, 1999.

\bibitem[BM04]{BrMr}
Jeffrey Brock and Dan Margalit.
\newblock Weil-{P}etersson isometries via the pants complex.
\newblock preprint, 2004.

\bibitem[Bro02]{Brkwpvs}
Jeffrey~F. Brock.
\newblock The {W}eil-{P}etersson visual sphere.
\newblock preprint, 2002.

\bibitem[Bro03]{Brkwp}
Jeffrey~F. Brock.
\newblock The {W}eil-{P}etersson metric and volumes of 3-dimensional hyperbolic
  convex cores.
\newblock {\em J. Amer. Math. Soc.}, 16(3):495--535 (electronic), 2003.

\bibitem[BY04]{Chn}
Chen Bo-Yong.
\newblock Equivalence of the {B}ergman and {T}eichm\"{u}ller metrics on
  {T}eichm\"{u}ller spaces.
\newblock preprint, 2004.

\bibitem[CY80]{ChYa}
Shiu~Yuen Cheng and Shing~Tung Yau.
\newblock On the existence of a complete {K}\"ahler metric on noncompact
  complex manifolds and the regularity of {F}efferman's equation.
\newblock {\em Comm. Pure Appl. Math.}, 33(4):507--544, 1980.

\bibitem[DKW00]{DW1}
Georgios Daskalopoulos, Ludmil Katzarkov, and Richard Wentworth.
\newblock Harmonic maps to {T}eichm\"uller space.
\newblock {\em Math. Res. Lett.}, 7(1):133--146, 2000.

\bibitem[DP86]{DPh}
Eric D'Hoker and D.~H. Phong.
\newblock On determinants of {L}aplacians on {R}iemann surfaces.
\newblock {\em Comm. Math. Phys.}, 104(4):537--545, 1986.

\bibitem[DW03]{DW2}
Georgios Daskalopoulos and Richard Wentworth.
\newblock Classification of {W}eil-{P}etersson isometries.
\newblock {\em Amer. J. Math.}, 125(4):941--975, 2003.

\bibitem[Gru01]{Gru}
Samuel Grushevsky.
\newblock An explicit upper bound for {W}eil-{P}etersson volumes of the moduli
  spaces of punctured {R}iemann surfaces.
\newblock {\em Math. Ann.}, 321(1):1--13, 2001.

\bibitem[Hej90]{Hejreg}
Dennis~A. Hejhal.
\newblock Regular {$b$}-groups, degenerating {R}iemann surfaces, and spectral
  theory.
\newblock {\em Mem. Amer. Math. Soc.}, 88(437):iv+138, 1990.

\bibitem[Hua04]{Zh2}
Zheng Huang.
\newblock On asymptotic {W}eil-{P}etersson geometry of {T}eichm\"{u}ller
  {S}pace of {R}iemann surfaces.
\newblock preprint, 2004.

\bibitem[IT92]{ImTan}
Y.~Imayoshi and M.~Taniguchi.
\newblock {\em An introduction to {T}eichm\"uller spaces}.
\newblock Springer-Verlag, Tokyo, 1992.
\newblock Translated and revised from the Japanese by the authors.

\bibitem[Iva88]{Ivcurv}
N.~V. Ivanov.
\newblock Teichm\"uller modular groups and arithmetic groups.
\newblock {\em Zap. Nauchn. Sem. Leningrad. Otdel. Mat. Inst. Steklov. (LOMI)},
  167(Issled. Topol. 6):95--110, 190--191, 1988.

\bibitem[Iva97]{Ivaut}
Nikolai~V. Ivanov.
\newblock Automorphisms of complexes of curves and of {T}eichm\"uller spaces.
\newblock In {\em Progress in knot theory and related topics}, volume~56 of
  {\em Travaux en Cours}, pages 113--120. Hermann, Paris, 1997.

\bibitem[KK04]{KKo}
Alexey Kokotov and Dmitry Korotkin.
\newblock Bergmann tau-function on {H}urwitz spaces and its applications.
\newblock preprint, 2004.

\bibitem[Kor99]{Krk}
Mustafa Korkmaz.
\newblock Automorphisms of complexes of curves on punctured spheres and on
  punctured tori.
\newblock {\em Topology Appl.}, 95(2):85--111, 1999.

\bibitem[KS93]{KShar}
Nicholas~J. Korevaar and Richard~M. Schoen.
\newblock Sobolev spaces and harmonic maps for metric space targets.
\newblock {\em Comm. Anal. Geom.}, 1(3-4):561--659, 1993.

\bibitem[LSY04]{LSY}
Kefeng Liu, Sun Sun, Xiaofeng, and Shing-Tung Yau.
\newblock Canonical metrics on the moduli space of {R}iemann surfaces, {I},
  {II}.
\newblock preprints, 2004.

\bibitem[Luo00]{Luaut}
Feng Luo.
\newblock Automorphisms of the complex of curves.
\newblock {\em Topology}, 39(2):283--298, 2000.

\bibitem[Mas76]{Msext}
Howard Masur.
\newblock Extension of the {W}eil-{P}etersson metric to the boundary of
  {T}eichmuller space.
\newblock {\em Duke Math. J.}, 43(3):623--635, 1976.

\bibitem[Mas01]{MasFN}
Bernard Maskit.
\newblock Matrices for {F}enchel-{N}ielsen coordinates.
\newblock {\em Ann. Acad. Sci. Fenn. Math.}, 26(2):267--304, 2001.

\bibitem[McM00]{McM}
Curtis~T. McMullen.
\newblock The moduli space of {R}iemann surfaces is {K}\"ahler hyperbolic.
\newblock {\em Ann. of Math. (2)}, 151(1):327--357, 2000.

\bibitem[McS98]{McSh}
Greg McShane.
\newblock Simple geodesics and a series constant over {T}eichmuller space.
\newblock {\em Invent. Math.}, 132(3):607--632, 1998.

\bibitem[Mir04a]{Mirgrow}
Maryam Mirzakhani.
\newblock Growth of the number of simple closed geodesics on hyperbolic
  surfaces.
\newblock preprint, 2004.

\bibitem[Mir04b]{Mirvol}
Maryam Mirzakhani.
\newblock Simple geodesics and {W}eil-{P}etersson volumes of moduli spaces of
  bordered {R}iemann surfaces.
\newblock preprint, 2004.

\bibitem[Mir04c]{Mirwitt}
Maryam Mirzakhani.
\newblock {W}eil-{P}etersson volumes and intersection theory on the moduli
  space of curves.
\newblock preprint, 2004.

\bibitem[MT04]{McTk}
Andrew McIntyre and Leon~A. Takhtajan.
\newblock Holomorphic factorization of determinants of {L}aplacians on
  {R}iemann surfaces and a higher genus generalization of {K}ronecker's first
  limit formula.
\newblock preprint, 2004.

\bibitem[MW02]{MW}
Howard Masur and Michael Wolf.
\newblock The {W}eil-{P}etersson isometry group.
\newblock {\em Geom. Dedicata}, 93:177--190, 2002.

\bibitem[MY83]{MoYa}
Ngaiming Mok and Shing-Tung Yau.
\newblock Completeness of the {K}\"ahler-{E}instein metric on bounded domains
  and the characterization of domains of holomorphy by curvature conditions.
\newblock In {\em The mathematical heritage of Henri Poincar\'e, Part 1
  (Bloomington, Ind., 1980)}, volume~39 of {\em Proc. Sympos. Pure Math.},
  pages 41--59. Amer. Math. Soc., Providence, RI, 1983.

\bibitem[Nag88]{Ng}
Subhashis Nag.
\newblock {\em The complex analytic theory of {T}eichm\"uller spaces}.
\newblock Canadian Mathematical Society Series of Monographs and Advanced
  Texts. John Wiley \& Sons Inc., New York, 1988.
\newblock A Wiley-Interscience Publication.

\bibitem[Obi99]{Obit1}
Kunio Obitsu.
\newblock Non-completeness of {Z}ograf-{T}akhtajan's {K}\"ahler metric for
  {T}eichm\"uller space of punctured {R}iemann surfaces.
\newblock {\em Comm. Math. Phys.}, 205(2):405--420, 1999.

\bibitem[Obi01]{Obit2}
Kunio Obitsu.
\newblock The asymptotic behavior of {E}isenstein series and a comparison of
  the {W}eil-{P}etersson and the {Z}ograf-{T}akhtajan metrics.
\newblock {\em Publ. Res. Inst. Math. Sci.}, 37(3):459--478, 2001.

\bibitem[OW]{WlO}
Kunio Obitsu and Scott~A. Wolpert.
\newblock in preparation.

\bibitem[Pen92]{Pen}
R.~C. Penner.
\newblock Weil-{P}etersson volumes.
\newblock {\em J. Differential Geom.}, 35(3):559--608, 1992.

\bibitem[Roy75]{Royicm}
H.~L. Royden.
\newblock Intrinsic metrics on {T}eichm\"uller space.
\newblock In {\em Proceedings of the International Congress of Mathematicians
  (Vancouver, B. C., 1974), Vol. 2}, pages 217--221. Canad. Math. Congress,
  Montreal, Que., 1975.

\bibitem[Sar87]{Sardet}
Peter Sarnak.
\newblock Determinants of {L}aplacians.
\newblock {\em Comm. Math. Phys.}, 110(1):113--120, 1987.

\bibitem[ST01]{SchTr}
Georg Schumacher and Stefano Trapani.
\newblock Estimates of {W}eil-{P}etersson volumes via effective divisors.
\newblock {\em Comm. Math. Phys.}, 222(1):1--7, 2001.

\bibitem[Tro86]{Trcurv}
A.~J. Tromba.
\newblock On a natural algebraic affine connection on the space of almost
  complex structures and the curvature of {T}eichm\"uller space with respect to
  its {W}eil-{P}etersson metric.
\newblock {\em Manuscripta Math.}, 56(4):475--497, 1986.

\bibitem[TT04a]{TakTI}
Lee-Peng Teo and Leon~A. Takhtajan.
\newblock {W}eil-{P}etersson metric on the universal {T}eichm\"uller space i:
  curvature properties and {C}hern forms.
\newblock preprint, 2004.

\bibitem[TT04b]{TakTII}
Lee-Peng Teo and Leon~A. Takhtajan.
\newblock {W}eil-{P}etersson metric on the universal {T}eichm\"uller space ii:
  {K}\"ahler potential and period mapping.
\newblock preprint, 2004.

\bibitem[TWZ04]{TanWZ}
Ser~Peow Tan, Yan~Loi Wong, and Ying Zhang.
\newblock Generalizations of {M}c{S}hane's identity to hyperbolic
  cone-surfaces.
\newblock preprint, 2004.

\bibitem[TZ91]{TkZg}
L.~A. Takhtajan and P.~G. Zograf.
\newblock A local index theorem for families of {$\overline\partial$}-operators
  on punctured {R}iemann surfaces and a new {K}\"ahler metric on their moduli
  spaces.
\newblock {\em Comm. Math. Phys.}, 137(2):399--426, 1991.

\bibitem[Wen01]{Weng}
Lin Weng.
\newblock {$\Omega$}-admissible theory. {II}. {D}eligne pairings over moduli
  spaces of punctured {R}iemann surfaces.
\newblock {\em Math. Ann.}, 320(2):239--283, 2001.

\bibitem[Wen04]{Wreg}
Richard~A. Wentworth.
\newblock Regularity of harmonic maps from riemann surfaces to the
  {W}eil-{P}etersson completion of {T}eich\"uller space.
\newblock preprint, 2004.

\bibitem[Wol]{Wlreprise}
Scott~A. Wolpert.
\newblock Convexity of geodesic-length functions: a reprise.
\newblock In {\em Spaces of Kleinian Groups}, Lec. Notes. Cambridge Univ.
  Press.
\newblock to appear.

\bibitem[Wol82]{WlFN}
Scott~A. Wolpert.
\newblock The {F}enchel-{N}ielsen deformation.
\newblock {\em Ann. of Math. (2)}, 115(3):501--528, 1982.

\bibitem[Wol85]{Wldtau}
Scott~A. Wolpert.
\newblock On the {W}eil-{P}etersson geometry of the moduli space of curves.
\newblock {\em Amer. J. Math.}, 107(4):969--997, 1985.

\bibitem[Wol86]{Wlchern}
Scott~A. Wolpert.
\newblock Chern forms and the {R}iemann tensor for the moduli space of curves.
\newblock {\em Invent. Math.}, 85(1):119--145, 1986.

\bibitem[Wol87a]{Wlselb}
Scott~A. Wolpert.
\newblock Asymptotics of the spectrum and the {S}elberg zeta function on the
  space of {R}iemann surfaces.
\newblock {\em Comm. Math. Phys.}, 112(2):283--315, 1987.

\bibitem[Wol87b]{Wlnielsen}
Scott~A. Wolpert.
\newblock Geodesic length functions and the {N}ielsen problem.
\newblock {\em J. Differential Geom.}, 25(2):275--296, 1987.

\bibitem[Wol90]{Wlhyp}
Scott~A. Wolpert.
\newblock The hyperbolic metric and the geometry of the universal curve.
\newblock {\em J. Differential Geom.}, 31(2):417--472, 1990.

\bibitem[Wol03]{Wlcomp}
Scott~A. Wolpert.
\newblock {G}eometry of the {W}eil-{P}etersson completion of {T}eichm\"{u}ller
  space.
\newblock In {\em Surveys in Differential Geometry VIII: Papers in Honor of
  Calabi, Lawson, Siu and Uhlenbeck}, pages 357--393. Intl. Press, Cambridge,
  MA, 2003.

\bibitem[Yam01]{Yam}
Sumio Yamada.
\newblock Weil-{P}etersson {C}ompletion of {T}eichm\"{u}ller {S}paces and
  {M}apping {C}lass {G}roup {A}ctions.
\newblock preprint, 2001.

\bibitem[Yeu03]{Yun1}
Sai-Kee Yeung.
\newblock Bounded smooth strictly plurisubharmonic exhaustion functions on
  {T}eichm\"uller spaces.
\newblock {\em Math. Res. Lett.}, 10(2-3):391--400, 2003.

\bibitem[Yeu04]{Yun}
Sai-Kee Yeung.
\newblock Quasi-isometry of metrics on {T}eichm\"{u}ller spaces.
\newblock {\em Int. Math. Res. Not.}, to appear, 2005.

\end{thebibliography}

\clearpage

\bigskip

\end{document}